\documentclass[reqno,oneside,12pt,hidelinks]{amsart}
\usepackage{hyphenat}
\usepackage{amssymb,amsmath,amsthm,amsxtra,calc,bm, color}
\usepackage{enumerate}
\usepackage[margin=.95in]{geometry}
\usepackage[scr=boondoxo]{mathalfa}
\usepackage{mathtools}

\usepackage{mleftright}
\mleftright

\usepackage{nccmath}
\usepackage{cases}

\usepackage{calc}
\usepackage{graphicx}
\usepackage{hyperref}
\usepackage[final]{microtype}

\newtheorem{theorem}{Theorem}[section]

\newtheorem{Property}[theorem]{Property}

\theoremstyle{remark}
\newtheorem*{remark}{Remark}

\newtheorem*{definition}{Definition}
\numberwithin{equation}{section}

\mathtoolsset{showonlyrefs}

\title{On Generalizations of the Cohen Structure Theorem}

\author{Robert Dicks}
\address{Department of Mathematics\\
Clemson University\\
Clemson, SC 29634} 
\email{rdicks@clemson.edu}

\begin{document}

\begin{abstract}
The Cohen structure theorem is an important tool in commutative algebra and algebraic geometry. Recent work of Caviglia, Goswami, and Mesiti addresses the related question of which unital rings admit a maximal ideal for which the natural projection admits a section; for such rings, they prove a structure theorem in terms of semidirect products of rings.

In this paper, we extract the underlying mechanism of their structure theorem to address when unital rings admit a nonzero proper ideal whose associated projection map has this property. We also address other settings involving ideals with conditions which are determined by ring theoretic properties of the corresponding quotient. 
\end{abstract}

\maketitle
\section{Introduction}
Let $R$ be a unital commutative complete Noetherian local ring with maximal ideal $\mathfrak{m}$. We say that the ring $R$ is \emph{equicharacteristic} if $R$ has the same characteristic as its residue field $R/\mathfrak{m}$. The Cohen structure theorem, which has numerous applications in commutative algebra and algebraic geometry (see e.g. \cite{Zariski} and \cite{Goto-Nishida}, along with the introduction of \cite{Caviglia-Goswami-Mesiti} for more discussion on uses of the Cohen structure theorem), states that $R$ has a subfield $F$ such that the image of $F$ under the natural projection map $R \rightarrow R/\mathfrak{m}$ is $R/\mathfrak{m}$; this property is equivalent to the natural projection map admitting a section.

Recent work of Caviglia, Goswami, and Mesiti 
\cite{Caviglia-Goswami-Mesiti} examines when unital (not necessarily commutative) rings have a maximal two-sided ideal whose associated projection map admits a section.
Among their main results is a structure theorem for such rings in terms of semidirect products of rings, \cite[Theorem $3.2$]{Caviglia-Goswami-Mesiti} (see Section $2$ for details on semidirect products of rings; these are the only semidirect products we consider) which was inspired by a communication between the second author and George Janelidze.
It is natural to consider whether the underlying mechanism of this result can be adapted to other settings. 

What is the importance of the condition of maximality in \cite[Theorem $3.2$]{Caviglia-Goswami-Mesiti}? The main point of this paper is that, in this discourse, maximality is not at all important;
We first pivot away from considering maximal ideals by examining unital rings $R$ which satisfy the following property.

\begin{Property}\label{mainproperty}
There exists a nonzero proper ideal I such that the natural projection map $R \rightarrow R/I$ admits a section.
\end{Property}
\begin{remark}
    All of the ideals that we consider are two-sided.
\end{remark}

 Our main result gives a structure theorem for unital rings which satisfy Property~\ref{mainproperty} in terms of semidirect products.

\begin{theorem}\label{maintheorem}
Let $R$ be a unital ring.
\begin{enumerate}
\item 
Suppose that $R=R_1 \rtimes R_2$ for rings $R_1$ and $R_2$ such that $R_1$ and $R_2$ are nonzero. Then $R$ satisfies Property~\ref{mainproperty} with an ideal $I$ such that $I \cong R_1$. 
\item
Suppose that $R$ satisfies Property~\ref{mainproperty} with an ideal $I$ and section $s:R/I \rightarrow R$. Then $R \cong I \rtimes s(R/I)$.  
\end{enumerate}

\end{theorem}
\begin{remark} 

(1)
Defining the semidirect product $I \rtimes s(R/I)$ requires specifying two structure maps $\lambda, \rho:s(R/I) \times I \rightarrow I$. As will be discussed in Section $3$, these maps will essentially be multiplication inside of $R$. \\
(2)
 The proof of the theorem shows that one may prove similar results when one places conditions on the ideal which may, roughly speaking, be determined by the ring theoretic properties of the corresponding quotient. For example, a proper ideal $I$ is a maximal left ideal if and only if $R/I$ is a division ring. As a porism to the proof of Theorem~\ref{maintheorem}, if we add the condition of $I$ being a maximal left ideal in Property~\ref{mainproperty}, then the result holds when we further assume in part $(1)$ that $R_2$ is a division ring; in part $(2)$, we can further conclude that $s(R/I)$ is a division ring. In other words, having Property~\ref{mainproperty} with the added condition of $I$ being a maximal left ideal is tantamount to saying that $R$ is a semidirect product of a nonzero ring and a division ring. If $R$ is commutative, then similar remarks apply when we further assume that $I$ is a prime ideal in Property~\ref{mainproperty} and that $R_2$ is an integral domain, or if we further assume that $I$ is a maximal ideal in Property~\ref{mainproperty} and that $R_2$ is a field. \\
(3)
 Each of the isomorphisms in Theorem~\ref{maintheorem} and its proof are ring isomorphisms.

\end{remark}

The paper is organized as follows. In the next section, we provide background on semidirect products of rings. In Section $3$ we prove our main result by  leveraging the central core of  the proof of \cite[Theorem $3.2$]{Caviglia-Goswami-Mesiti} in this domain. The proof of part $(1)$ leverages the projection map onto the second factor; its kernel is an ideal which witnesses the fact that Property~\ref{mainproperty} holds. The proof of part $(2)$ proceeds by considering the natural map which identifies the abelian group $I \oplus s(R/I)$ with a subgroup of $R$ (the map $\Phi:I \rtimes s(R/I) \rightarrow R$ below). The semidirect product structure on $I \oplus s(R/I)$ upgrades $\Phi$ to a ring homomorphism. The fact that $s$ is a section allows one to produce an inverse $\Psi:R\rightarrow I \rtimes s(R/I)$ to $\Phi$ by splitting $r \in R$ into its ``$I$-component" $r-s([r])$ and its ``$s(R/I)$-component" $s([r])$.
We also discuss how to prove this result in the other settings mentioned in part $(2)$ of the remark following Theorem~\ref{maintheorem}.  
\section{Background/Conventions}
Throughout the paper, we assume all the rings which we consider are associative but not necessarily unital. If $R$ is a ring, we write $R_R$ (respective $_RR$) when we consider $R$ as a right $R$-module (respectively left $R$-module) over itself in the natural way. For an $R$-module $M$ (which may be left or right), we write $\text{End}_{\text{$R$-Mod}}(M)$ for the ring of $R$-module endomorphisms of $M$. 

We recall the notion of a semidirect product of rings.
\begin{definition}
Let $R$ and $T$ be rings. Let $\lambda:R\rightarrow \text{End}_{\text{$T$-Mod}}(T_T)$ be a ring homomorphism and $\rho:R\rightarrow \text{End}_{\text{$T$-Mod}}(_TT)$ be a ring anti-homomorphism such that, for every $s,t \in R$ and $x,y \in T$, we have
\begin{equation}\label{compatibility1}
\lambda(s) \circ \rho(t)=\rho(t) \circ \lambda(s) 
\end{equation}
\begin{equation}\label{compatibility2}
\rho(s)(x)y=x\lambda(s)(y).
\end{equation}
The \emph{semidirect product} $T \rtimes_{\lambda, \rho}R$ with respect to $\lambda$ and $\rho$ is defined as follows. The underlying abelian group is given by $T \oplus R$ and we define the product of $(t_1,r_1),(t_2,r_2) \in T \oplus R$ via
\[
(t_1,r_1) \cdot (t_2,r_2):=(t_1t_2+\lambda(r_1)(t_2)+\rho(r_2)(t_1), r_1r_2).
\]
\end{definition}
\begin{remark}
    It is straightforward to check that the semidirect product is a ring. We also have ring homomorphisms $i_T:T \rightarrow T \rtimes_{\lambda, \rho} R$ and $i_R:R \rightarrow T \rtimes_{\lambda, \rho}R$ given by $i_T(t)=(t,0)$ and $i_R(r)=(0,r)$.
\end{remark}
If $R$ is a ring and $I \subseteq R$ is an ideal, we denote the image of $r \in R$ under the natural projection map $R \rightarrow R/I$ by $[r]$. 
\section{Proof of Theorem~\ref{maintheorem}/generalizations}
\subsection{Proof of Theorem~\ref{maintheorem}}
We begin the proof of part $(1)$. Consider the map 
\[
\phi:R_1 \rtimes_{\lambda, \rho}R_2\rightarrow R_2
\]
given by $\phi(r_1,r_2)=r_2$. We have
\[
\text{ker}(\phi)=\{(r_1,0) \ : \ r_1 \in R_1\} \cong R_1;
\]
henceforth, we identify $\text{ker}(\phi)$ with $R_1$. With this notation, we have
\begin{equation}\label{refinementfuel}
R_1 \rtimes_{\lambda,\rho}R_2/R_{1}\cong R_2.
\end{equation}
Composing this map with $i_{R_2}$ provides a section $s:R_1 \rtimes_{\lambda,\rho} R_2/R_1 \rightarrow R_1 \rtimes_{\lambda,\rho} R_2$ to the natural projection map $R_1 \rtimes_{\lambda,\rho} R_2 \rightarrow R_1 \rtimes_{\lambda,\rho}R_2/R_1$. Indeed, we have $s([(x,u)])=(0,u)$, and $[(0,u)]=[(x,u)]$, so part $(1)$ holds.

We now prove part $(2)$. It is simple to check that the map $\lambda:s(R/I) \rightarrow \text{End}_{\text{$R$-Mod}}(I_I)$ given by $\lambda(x)(b)=xb$ is a ring homomorphism and that the map $\rho:s(R/I) \rightarrow \text{End}_{\text{$R$-Mod}}(_II)$  given by $\rho(x)=bx$ is a ring anti-homomorphism (and also that, for each $x \in s(R/I)$, the maps $\lambda(x):I \rightarrow I$ and $\rho(x):I \rightarrow I$ are $R$-module endomorphisms) satisfying the conditions \eqref{compatibility1} and \eqref{compatibility2}; henceforth, we suppress the subscripts by writing $I \rtimes s(R/I)$.

We show that $I \rtimes s(R/I) \cong R$. Consider the map
\[
\Phi:I \rtimes s(R/I) \rightarrow R
\]
given by $\Phi(r_1,r_2)=r_1+r_2$. It is not hard to check that $\Phi$ is a ring homomorphism.

Consider now the map $\Psi:R \rightarrow I \rtimes s(R/I)$ given by $\Psi(r)=(r-s([r]),s([r]))$ (since $[r-s([r])]=[r]-[s([r])]=[r]-[r]=[0]$, we have $r-s([r]) \in I$). It is straightforward to check that $\Psi$ is a ring homomorphism.

We now show that $\Phi$ and $\Psi$ are inverses to each other. We have
\[
\Phi(\Psi(r))=\Phi(r-s([r]),s([r]))=r.
\]
To show that $\Psi \circ \Phi=\text{id}$, observe that 
\[
\Psi(\Phi(r_1,r_2))=(r_1+r_2-s([r_1+r_2]),s([r_1+r_2])),
\]
so we need only show that $s([r_1+r_2])=r_2$. Observe that $r_2=s([r_2])$, and since $s$ is injective by virtue of being a section, it now suffices to show that $[r_1+r_2]=[r_2]$, which is clear since $r_1 \in I$. The result follows. 
\subsection{Other settings}
We can refine the above proof in several ways. What follows is largely the same argument as in Section $3.1$; we sketch the details and keep the same notation. 

Recall that a proper ideal in a unital ring is a maximal left ideal if and only if the corresponding quotient is a division ring. In part $(1)$ of Theorem~\ref{maintheorem}, if we further assume that $R_2$ is a division ring, then $\text{ker}(\phi) \cong R_1$ is a left maximal ideal in $R$ for which the natural projection map has a section via the map in \eqref{refinementfuel}. If we assume in part $(2)$ that $I$ is a maximal left ideal, then the isomorphism argument still works as above, but we further see that $s(R/I)$ is a division ring.

If we assume that $R$ is commutative, then the same argument applies when we add the condition in Property~\ref{mainproperty} that $I$ is prime and assume that $R_2$ is an integral domain.
\section{AI statement}
The author did not use generative AI in the production of his arguments nor in their exposition. The author used ChatGPT to proofread the paper.
\bibliography{Cohenstructure}

@article{Goto-Nishida,
	author = {Goto, Shiro and Nishida, Kenji},
	doi = {10.1017/S0024610700001927},
	fjournal = {Journal of the London Mathematical Society. Second Series},
	issn = {0024-6107,1469-7750},
	journal = {J. London Math. Soc. (2)},
	mrclass = {13D05},
	mrnumber = {1810132},
	mrreviewer = {Carl\ Faith},
	number = {2},
	pages = {319--335},
	title = {Finite modules of finite injective dimension over a {N}oetherian algebra},
	url = {https://doi.org/10.1017/S0024610700001927},
	volume = {63},
	year = {2001}}

@article{Zariski,
	author = {Zariski, Oscar},
	doi = {10.2307/1990628},
	fjournal = {Transactions of the American Mathematical Society},
	issn = {0002-9947,1088-6850},
	journal = {Trans. Amer. Math. Soc.},
	mrclass = {14.0X},
	mrnumber = {21694},
	mrreviewer = {Andr\'{e}\ Weil},
	pages = {1--52},
	title = {The concept of a simple point of an abstract algebraic variety},
	url = {https://doi.org/10.2307/1990628},
	volume = {62},
	year = {1947}}

@article{Caviglia-Goswami-Mesiti,
	author = {Caviglia, Elena and Goswami, Amartya and Mesiti, Luca},
	doi = {10.1016/j.jalgebra.2026.05.014},
	fjournal = {Journal of Algebra},
	issn = {0021-8693,1090-266X},
	journal = {J. Algebra},
	mrclass = {13H10 (13B30 13E05)},
	mrnumber = {5083163},
	pages = {293--310},
	title = {On extensions of {C}ohen structure theorem},
	url = {https://doi.org/10.1016/j.jalgebra.2026.05.014},
	volume = {707},
	year = {2026}}
\bibliographystyle{amsalpha}
\end{document}